\documentclass[a4paper, 11pt]{amsart}

\usepackage[T1]{fontenc}
\usepackage[utf8]{inputenc}

\usepackage{amssymb}

\makeatletter
\@namedef{subjclassname@2020}{%
  \textup{2020} Mathematics Subject Classification}
\makeatother

\allowdisplaybreaks[2]

\theoremstyle{definition}

\begin{document}

\title[Frobenius' Theorem]{Frobenius' Theorem on  Division Algebras}
\author{Matej Bre\v sar}

\address{Faculty of Mathematics and Physics,  University of Ljubljana $\&$
  Faculty of Natural Sciences and Mathematics, University 
of Maribor $\&$ IMFM, Ljubljana, Slovenia}
 \email{matej.bresar@fmf.uni-lj.si}

\thanks{Partially supported by ARIS Grant P1-0288 }

\keywords{Frobenius' Theorem, quaternions, division algebra}

\subjclass[2020]{16K20}

\begin{abstract} Frobenius' Theorem states that the only finite-dimensional real division algebras are the algebra of real numbers $\mathbb R$, the algebra of complex numbers $\mathbb C$, and the algebra of quaternions $\mathbb H$. We present a short proof which uses only standard undergraduate mathematics.
\end{abstract} 
\maketitle

Recall that the {\bf algebra of quaternions} $\mathbb H$
 is a 
$4$-dimensional real  algebra whose standard basis consists of unity $1$ and elements $i,j,k$ satisfying
\begin{equation}\begin{aligned} \label{h}&i^2=j^2=k^2=-1,\\
ij=-ji=k,&\,\,\, jk=-kj =i,\,\,\,ki=-ik=j.\end{aligned}\end{equation}
These identities imply that
\begin{equation}\label{n}
(\lambda_0 + \lambda_1 i + \lambda_2 j + \lambda_3 k)(\lambda_0 - \lambda_1 i - \lambda_2 j - \lambda_3 k) = \lambda_0^2+
\lambda_1^2+\lambda_2^2+\lambda_3^2\end{equation}
for all $\lambda_i\in\mathbb R$,
from which it readily follows that $\mathbb H$ is a division algebra. It is by far the most well known
 noncommutative division algebra, 
occurring  throughout mathematics and its applications.
The notation $\mathbb H$ is in honor of  William Rowan Hamilton, who discovered the quaternions in 1843.

In 1878,  Ferdinand Georg  Frobenius proved that $\mathbb{R}$, $\mathbb C$, and $\mathbb H$ are the only  finite-dimensional real division algebras \cite{F}. In modern literature, this theorem is often stated for 
 (seemingly more general) {\bf algebraic 
real division algebras}; we say that a real algebra $D$ is algebraic if, for every $x\in D$, there exists a nonzero polynomial
 $f(X)\in \mathbb R[X]$  such that $f(x)=0$. Every finite-dimensional algebra $D$ is algebraic. Indeed, if $D$
is $n$-dimensional, then the elements $1,x,\dots,x^n$ are linearly dependent, meaning that $x$ is a root of 
 a nonzero real  polynomial of degree at most $n$.

A number of different proofs of Frobenius' Theorem, some elementary and some advanced, appear in the literature. We refer to \cite{INCA}, \cite{BS}, and \cite{P} for three
 different short and elementary proofs.
Combining  ideas from \cite{INCA} and \cite{BS}, we will give yet another such proof, which is somewhat 
 more conceptual than those from \cite{INCA} and \cite{BS}, and  therefore easier to memorize. Also, it may  be interesting for students to see results  they know from  other mathematical areas used in this proof.
The following  will be  our main tools:
\begin{itemize}
\item[(a)] Every nonconstant polynomial $f(X)\in \mathbb R[X]$ 
is the  product of linear and quadratic polynomials in $\mathbb R[X]$ (this is a version of the Fundamental Theorem of Algebra).
\item[(b)] Every linear operator from an odd-dimensional real vector space to itself has a real eigenvalue.
\item[(c)] If $U$ is a proper subspace of a finite-dimensional real inner product vector space $V$, then $V$ contains a unit vector $e$ orthogonal to all vectors in $U$.
\end{itemize}

Before proceeding to the proof, we 
record some notational, along with  motivational remarks.
Let $D$ be any real division algebra.  We write $\lambda 1\in D$, where $\lambda$ is a real number, simply as $\lambda$, and accordingly identify $\mathbb R1$ with $\mathbb R$. Thus,
by saying that $x\in D$ satisfies $x\le 0$ we mean that $x=\lambda 1$ with $\lambda$  a nonpositive real number.
 We write 
$$x\circ y= xy+yx$$
for  any   $x,y\in D$. 
 Observe that  $x\circ y=y\circ x$, $x\circ (y+z)=x\circ y + x\circ z$, and
 \begin{equation} \label{e} x\circ y = (x+y)^2-x^2-y^2.\end{equation}
If $x\circ y=0$, then we say that $x$ and $y$ {\bf anticommute}.

 We remark that the quaternions $i$, $j$, and $k$ anticommute with each other. Note also that
$V= {\rm span}\{i,j,k\}$, the linear span of $i,j,k$, is a linear subspace of $\mathbb H$ such that $\mathbb H =\mathbb R\oplus V$, and any of its  elements
$u=\lambda_1 i + \lambda_2 j + \lambda_3 k$ and $v=\mu_1 i + \mu_2 j + \mu_3 k$
 satisfy  
$$ -\frac{1}{2}u\circ v= \lambda_1\mu_1 + \lambda_2\mu_2 + \lambda_3\mu_3\in\mathbb R\quad\mbox{and}\quad v^2 = -(\mu_1^2 + \mu_2^2 + \mu_3 ^2)\le 0.$$
These observations provide the intuition behind our proof, which we  now  give.


\bigskip \noindent
 {\bf Frobenius' Theorem.} {\em  An algebraic  real division algebra $D$ is isomorphic  to 
$\mathbb R$, $\mathbb C$, or $\mathbb H$.}
\smallskip

\begin{proof} Set  $V=\{v\in D\,|\,v^2\le 0\}$. 
We claim that for any $x\in D$,
\begin{equation}\label{a} x^2\in \mathbb R\implies x\in \mathbb R \mbox{ or } x\in V.\end{equation}
Indeed, if $x\notin V$, then $x^2 >0$ and hence $x^2= \gamma^2$ with $\gamma\in\mathbb R$, which gives  $(x-\gamma)(x + \gamma)=0$ and so
 $x=\pm\gamma\in\mathbb R$.

Next, we claim  that for every $x\in D$ there exists an
$\alpha\in\mathbb R$ such that
\begin{equation} \label{rv} x^2 +\alpha x \in\mathbb R \quad\mbox{and}\quad
x + \frac{\alpha}{2}\in V.
\end{equation}
We may assume that $x\notin \mathbb R$. Let $f(X)\in \mathbb R[X]$ be a nonconstant polynomial such that $f(x)=0$. Write $f(X)$ as in (a). As $D$ has no zero divisors and $x\notin\mathbb R$, it follows that
$x^2 +\alpha x + \beta =0$ for some $\alpha,\beta\in\mathbb R$. Hence,  
$(x + \frac{\alpha}{2})^2 \in \mathbb R$. From \eqref{a} we conclude that $x + \frac{\alpha}{2}\in V$, which completes the proof of \eqref{rv}.  

Our next goal is to prove that 
\begin{equation}\label{uv} u\circ v \in \mathbb R\,\,\mbox{ for all $u,v\in V$.}
\end{equation}
We may assume that $u$ and $v$ are linearly independent.
By squaring both sides we see that 
$v=\lambda +\mu u$ with $\lambda,\mu\in \mathbb R$  leads to a contradiction, so $W=
{\rm span}\{1,u,v\}$ is a $3$-dimensional space. Further,  \eqref{rv}
shows that $x^2\in W$ for all $x\in W$,
and hence, by \eqref{e}, 
 $x\circ y\in W$ for all $x,y\in W$. Consequently, the 
 linear operator $$x\mapsto x\circ v$$ maps  $W$ to itself. 
By (b), there exist a $\lambda\in\mathbb R$ and a  $w\ne 0$ in $W$ such that
$w\circ v = \lambda w$. That is, $wvw^{-1}= \lambda -v$. Squaring both sides and using $v^2\in\mathbb R$ 
we obtain $-2\lambda v\in \mathbb R$ and hence
$\lambda=0$ since $v\notin \mathbb R$.
 Thus, $w\circ v=0$.   Hence $w$ and $v$ do not commute, so writing 
$w$ as $\alpha_0 + \alpha_1 u + \alpha_2 v$ we see that $\alpha_1\ne 0$, from which it follows that 
$u\circ v\in {\rm span}\{1,v\}$. As $u$ and $v$ occur symmetrically,  we analogously have $u\circ v\in {\rm span}\{1,u\}$. The linear independence of $1,u,v$ implies \eqref{uv}.

We continue to assume that $u$ and $v$ are linearly independent elements in $V$.
In light of \eqref{uv},  \eqref{e}  shows that  $(u+v)^2\in\mathbb R$. Since, as shown in the preceding paragraph, $1,u,v$ are linearly independent,  
\eqref{a} tells us that $u+v\in V$. As $V$ is clearly closed under scalar multiplication, this proves that 
it is a linear subspace of $D$.  By \eqref{rv}, $D = \mathbb R\oplus V$.

Observe that \eqref{uv} implies that
$$\langle u\,|\,v\rangle = -\frac{1}{2}u\circ v$$ defines an inner product on $V$.
We may assume that 
$V\ne \{0\}$. Let $i$ be a unit vector in $V$. Then $i^2=-1$, so  $D\cong\mathbb C$ if $\dim_\mathbb R V = 1$. Let
$\dim_\mathbb R V > 1$. By (c), there is a unit vector $j\in V$ orthogonal to $i$. Set $k=ij$. One immediately checks that $i,j,k$ satisfy \eqref{h}. Moreover,
since \eqref{h} implies \eqref{n}, 
$1,i,j,k$ are linearly independent. Thus, $D\cong \mathbb H$ if  
$\dim_\mathbb R V =3$. Suppose 
$\dim_\mathbb R V >3$. Using (c) again, we see that there is a unit vector $e\in V$ orthogonal to span$\{i,j,k\}$. This means that 
each of $i$, $j$, $k$
anticommutes with $e$. 
However, since $k=ij$ and the product of two elements anticommuting with $e$ obviously commutes with $e$, we arrive at a contradiction that
$ek=0$. 
\end{proof}

An obvious corollary of Frobenius' Theorem is that an odd-dimensional real division algebra $D$ is isomorphic to $\mathbb R$. As noticed in \cite{BS}, this  immediately follows from (b). Indeed, 
given any $d\in D$, the linear operator $x\mapsto dx$ must have a real eigenvalue $\lambda$, hence $(d-\lambda)w =0$ for some  $0 \ne w\in D$, and so $d=\lambda\in\mathbb R$. This little observation may be useful for motivating the introduction of the quaternions. (Note that the associativity of multiplication was not used in this argument;
it is a classical result that
 nonassociative finite-dimensional real division algebras exist only in dimensions $1$, $2$, $4$, and $8$ \cite{BM, Ker}).

\end{document}